\documentclass{ifacconf}
\usepackage{amsmath,amssymb,amsfonts,mathtools}
\usepackage{natbib}
\usepackage{xcolor}
\usepackage{bm}
\usepackage{listings} 
\usepackage{graphicx}  
\usepackage{subcaption}
\newtheorem{remark}{Remark}
\newtheorem{example}{Example}

\newtheorem{theorem}{Theorem}

\begin{document}

\begin{frontmatter}

\title{Towards Programming Adaptive Linear Neural Networks Through Chemical Reaction Networks
\thanksref{footnoteinfo}} 

\thanks[footnoteinfo]{This work was funded by the National Nature Science Foundation of China under Grant No. 12071428 and 62111530247, and the Zhejiang Provincial Natural Science Foundation of China under Grant No. LZ20A010002.}

\author[First]{Yuzhen Fan,} 
\author[First]{Xiaoyu Zhang,} 
\author[First]{Chuanhou Gao} 

\address[First]{School of Mathematical Sciences, Zhejiang
University, Hangzhou 310027, China(e-mail: \{12035034,Xiaoyu\_Z,gaochou\}@zju.edu.cn)..}

\begin{abstract}
This paper is concerned with programming adaptive linear neural networks (ALNNs) using chemical reaction networks (CRNs) equipped with mass-action kinetics. Through individually programming the forward propagation and the backpropagation of ALNNs, and also utilizing the permeation walls technique, we construct a powerful CRN possessing the function of ALNNs, especially having the function of automatic computation. We also provide theoretical analysis and a case study to support our construction. The results will have potential implications for the developments of synthetic biology, molecular computer and artificial intelligence.

\end{abstract}


\end{frontmatter}

\section{Introduction}
Living cells can process and respond to complex information from the environment because of the sophisticated behavior resulting from chemical reactions between molecules. Recently, there is an increasing interest in designing, programming, and assembling biomolecular circuits to accomplish specific tasks in living cells in synthetic biology. Some gene expression circuits can compute logarithm, and cell receptor allosteric protein acts as a logarithm sensor in the cellular signaling pathway \citep{chou2017chemical}. Moreover, biological circuit design based on control theory is also a way to control the behavior of cells \citep{del2016control}, which could be used in diseases diagnosis and drug delivery. 

Chemical reaction networks (CRNs) are widely used to model biomolecular interactions, and arbitrary CRN with mass-action kinetics is demonstrated to be implemented as DNA strand displacement \citep{soloveichik2010dna}. It means that CRNs can be regarded as a bridge to design biological circuits. However, it is challenging to compile CRNs directly to accomplish complex tasks, such as pattern recognition, classification, and decision-making. Artificial neural networks provide a mature mathematical framework to address this problem in-silico. Therefore, we aim to program neural networks through biochemical reaction networks, which can be used as an assembled module in-vivo. There has much significant work in this study. \cite{hjelmfelt1991chemical} confirmed that a class of metabolic reactions behave like the firing and quiescent states in McCulloch-Pitts neuron.
The perceptron model with ReLu or Sigmoid activation functions was implemented by the molecular sequestration reaction and phosphorylation-dephosphorylation cycles, respectively \citep{Moorman2019A, samaniego2021signaling}. The implementation of the Binary-Weight ReLu neural networks by rate-independent CRNs was proposed in \citep{vasic2020deep}. Additionally, \cite{anderson2021reaction} provided a mathematical framework for the construction of deterministic reaction networks that implement neural networks. However, the design of CRNs above is for a neural network with given parameters, and there are fewer results on the neural network implementation capable of learning \citep{blount2017feedforward, chiang2015reconfigurable, banda2014training}, in which reactions were designed qualitatively. In addition, in the previous design, human intervention is inevitable, which results in the failure of these reaction systems to achieve automatic computation.
For these reasons, this paper tries to program a complete single-layer adaptive linear neural network (ALNN) using CRNs. The biochemical adaptive linear network model is capable of learning quantitatively as an implementation of an adaptive linear element, which achieves both forward propagation and the update of weights. The steady states of our reaction networks are defined as the computation results. Our work contributes to the embedding of adaptive computation in molecular contexts. The remainder of this paper is organized as follows. Section 2 gives preliminaries on CRNs and ALNNs. In Section 3, we use CRNs to program ALNNs, and the corresponding theoretical supports are also provided. Finally, Section 4 illustrates our methods through a linear fitting problem and further presents some discussions about the challenges, significance, and points of future study of the current work.

\section{Preliminaries}
In this section, we briefly introduce CRNs \citep{horn1972general,feinberg1972complex} and ALNNs \citep{widrow1960adaptive}.

\subsection{CRNs}
 A CRN consists of a finite set of \textit{species} $\mathcal{S}=\{X_{1}, \ldots,X_{n}\}$, a set of \textit{complexes} $\mathcal{C}=\{C_{1},\ldots,C_{c}\}$ and a set of \textit{reactions} $\mathcal{R}=\{R_{1},\ldots,R_{r}\}$ describing the interactions between species, often referred to as a triple ($\mathcal{S,C,R}$). The $j$th ($j=1,\ldots,r$) reaction is expressed as  
\begin{equation}
R_j:~~~~~	\sum_{i=1}^{n} \alpha_{ij} X_{i} \longrightarrow \sum_{i=1}^{n} \beta_{ij} X_{i},
	\label{eq:1}
\end{equation}
where the linear combination of species on the two sides of the arrow are called \textit{reactant complex} and \textit{product complex}, respectively, and the nonnegative integers $\alpha_{ij}, \beta_{ij}$ are the stoichiometric coefficients of $X_i$. The dynamics of CRNs usually focuses on describing the change of the concentrations of every species as time passes, denoted by $\tilde{x}_i$. When a CRN is assigned a mass-action kinetics, the rate for the reaction (\ref{eq:1}) is evaluated by $k_j\prod_{i=1}^{n} \tilde{x}_{i}^{\alpha_{ij}}$, where $k_j$ is the reaction rate constant of $R_j$. The dynamic equation of the CRN ($\mathcal{S,C,R}$) is thus expressed as 
\begin{equation}
	\frac{d\tilde{x}_{i}(t)}{dt} = \sum_{j=1}^{r} (\beta_{ij}-\alpha_{ij})k_j \prod_{l=1}^{n} \tilde{x}_{l}^{\alpha_{lj}},
	\label{eq:2} 
\end{equation}
where $i=1,2,\ldots,n$. Essentially, these ordinary differential equations (ODEs) offer a polynomial system. We call the quaternity  ($\mathcal{S,C,R},k$) a mass-action system (MAS). 

Turing universality asserts that any computation can be embedded into a class of polynomial ODEs \citep{fages2017strong}, and then these ODEs can be implemented with mass-action chemical kinetics. Note that the reactions to perform a certain computation are usually artificial, but they can be made to have biological significance if the DNA strand displacement technique \citep{soloveichik2010dna} is applied. This together with the strong computing power of CRNs encourages us to use them to program machine learning algorithms. We select the most fundamental one, ALNN, as a trial.


\subsection{ALNNs}
ALNNs were first proposed by Widrow and Hoff (1960), whose main purpose is to linearly approximate a function relationship. Structurally, an ALNN comprises an input layer with several nodes and an output layer with a node. The connection of the former to the latter is modeled by weights. Finally, a linear activation function controls the output. Mathematically, suppose that $\bm{x}=(x_{1},\ldots,x_{n})^{\top} \in \mathbb{R}^{n}$ is an input vector, and $\bm{\omega}=(\omega_1,\ldots,\omega_n)^{\top} \in \mathbb{R}^{n}$ is a weight vector, then the output (here we only consider the case of continuous analog output) is 
\begin{equation}
	y=\sum_{i=0}^{n} x_{i} \omega_{i}.
	\label{eq:3}
\end{equation}
Here, $\omega_0$ is the threshold of the output node, which could be regarded as a dumb node with input $x_{0}=1$.   

The learning ability of ALNN is reflected in that the weights can be updated. The frequently-used update algorithm is the least mean square (LMS) that changes the weights by minimizing the square deviation between the desired output $d_j$ and the real output $y_j$, written as $\frac{1}{2}(d_j-y_j)^2$, for every sample, i.e., the change formula following
\begin{equation}
	\Delta \omega_{ij}(m) = \eta x_{ij} (d_{j}-y_{j}(m)),
	\label{eq:4}
\end{equation}
where $i$ denotes the $i$th node in the input layer, $j$ represents the $j$th sample, $m$ indicates the $m$th iteration, and $\eta \in \left(0,1\right]$ is the learning rate.

The main procedures of ALNN are divided into assignment, feedforward, learning, judgment, and loop, and the involved computation includes addition, subtraction, and multiplication. The former two procedures belong to the forward propagation process while the latter three procedures belong to the backpropagation process. We use the CRN language to program these two processes subject to minimizing human intervention.

\section{Programming ALNNs}
In this section, we will implement ALNNs through CRNs according to the forward and backpropagation processes.  

\subsection{Programming Forward Propagation}
The key procedure in the forward propagation process of ALNNs is to compute the output based on the inputs and the weights, i.e., Eq. (\ref{eq:3}). To program it using CRNs, we need to correspond each variable in this equation to a species concentration of the designed CRN. Note a fact that the inputs, weights, and output can be any real values while the species concentration should be nonnegative. To make the right match, we use the dual-rail encoding technique \citep{vasic2020deep}, through which two species concentrations will correspond to a variable in ALNNs. Following this rule, Eq. (\ref{eq:3}) thus changes to be
\begin{equation}
	\begin{split}
		y &=\sum_{i=0}^{n} (x_{i}^{+}-x_{i}^{-})(\omega_{i}^{+}-\omega_{i}^{-})\\
		&=\sum_{i=0}^{n} (x_{i}^{+}\omega_{i}^{+}+x_{i}^{-}\omega_{i}^{-})-\sum_{i=0}^{n}(x_{i}^{+}\omega_{i}^{-}+x_{i}^{-}\omega_{i}^{+})\\
		&=y^{+}-y^{-},
	\end{split}
\label{eq:5}
\end{equation} 
where $y^+,y^-$ represent the first and the second part of the second equality, respectively, and $x_{0}^{+}-x_{0}^{-}=1$. The current pattern can realize the right match between the nonnegative species concentrations and the variables in ALNNs, e.g., $x_i^+,x_i^-$ are both nonnegative, but their difference $x_i^+-x_i^-$ can have any sign, and can thus match $x_i$. Physically speaking, the variables with ``+" superscript represent increment while those with ``-" superscript mean decrement, which also applies to all of the subsequent variables.   

Eq. (\ref{eq:5}) suggests that we need to introduce the \textit{input species pair} $\{X_{i}^{+}, X_{i}^{-}\}$, the \textit{weight species pair} $\{W_{i}^{+}, W_{i}^{-}\}$, and the \textit{output species pair} $\{Y^{+}, Y^{-}\}$ to build the CRN. The expressions of $y^+$ and $y^-$ indicate the equilibrium moulds of $Y^+$ and $Y^-$, i.e., values evaluated at $t=\infty$, which suggests a group of catalytic reactions needed to be constructed. The balancing of $Y^+$ and $Y^-$ means it necessary to build decay reactions in the network. Finally, the automatic computation of the output requests that the output can be only represented by the equilibrium concentration of $Y^+$ or the one of $Y^-$. This implies an annihilation reaction about these two species should be included. So, we construct the following MAS to program the forward propagation process,  
\begin{align}
	W_{i}^{+}+X_{i}^{+} & \stackrel{k}{\longrightarrow} W_{i}^{+}+X_{i}^{+}+ Y^{+}, \label{eq:6}\\
	W_{i}^{-}+X_{i}^{-} & \stackrel{k}{\longrightarrow} W_{i}^{-}+X_{i}^{-}+ Y^{+}, \label{eq:7}\\
	W_{i}^{+}+X_{i}^{-} & \stackrel{k}{\longrightarrow} W_{i}^{+}+X_{i}^{-}+ Y^{-}, \label{eq:8}\\
	W_{i}^{-}+X_{i}^{+} & \stackrel{k}{\longrightarrow} W_{i}^{-}+X_{i}^{+}+ Y^{-}, \label{eq:9}\\
	Y^{+} & \stackrel{k}{\longrightarrow} \varnothing, \label{eq:10}\\
	Y^{-} & \stackrel{k}{\longrightarrow} \varnothing, \label{eq:11}\\
	Y^{+}+Y^{-} & \stackrel{k_{\infty}}{\longrightarrow} \label{eq:12} \varnothing,
\end{align}
where $i=1,2,\ldots,n$, the first four reactions are catalytic reactions, the reactions (\ref{eq:10}), (\ref{eq:11}) are decay reactions while the last one is the annihilation reaction. Regarding the reaction rate constants $k$ and $k_{\infty}$, we set $k=1$ for every reaction for simplicity while $k_{\infty}$ is set such that $k_{\infty} \gg k$. The main reason is that we need the species $Y^+$ and $Y^-$ (produced by other reactions) to consume each other instantly in the rapid annihilation reaction (\ref{eq:12}). The final result is that a certain species is completely consumed and cannot be regenerated, and the other one has left whose equilibrium concentration is directly the output $y$.  

It should be noted that Moorman, Samaniego, et al. (2019) also designed a group of reactions for the perceptron model for the same purpose. Compared to theirs, the current first four catalytic reactions are different from their constitutive activation reactions while others are the same. However, the current design can perform an automatic calculation of Eq. (\ref{eq:5}) only based on the initial values of input species and weight species, whereas in those reactions \citep{Moorman2019A} it needs to manually calculate the sum of the product of inputs and weights. Moreover, ours can handle the case of negative output, which is also not covered by the work of Moorman, Samaniego, et al. (2019).

In the MAS given by the reactions (\ref{eq:6})-(\ref{eq:12}), the species $W_{i}^{+}, W_{i}^{-},X_{i}^{+}$ and $X_{i}^{-}$ are all catalysts, whose concentrations will keep unchanged as time passes, i.e., preserved at $\omega_i^+(0),\omega_i^-(0),x_i^+(0)$ and $x_i^-(0)$, respectively, so the dynamics of this MAS follows
\begin{equation}  
    \begin{split}
    	   \frac{dy^{+}(t)}{dt} & = kp^{+}(t)-ky^{+}(t)-k_{\infty}y^{+}(t)y^{-}(t), \\
    	   \frac{dy^{-}(t)}{dt} &= kp^{-}(t)-ky^{-}(t)-k_{\infty}y^{+}(t)y^{-}(t),
    	   \end{split}
    	   \label{eq:13}
\end{equation}
where
\begin{gather}\label{eq:00}
\begin{split}
	p^{+}(t)=\sum_{i=0}^{n}
	(x_{i}^{+}(t)\omega_{i}^{+}(t)+x_{i}^{-}(t)\omega_{i}^{-}(t)), \\
	p^{-}(t)=\sum_{i=0}^{n}
	(x_{i}^{+}(t)\omega_{i}^{-}(t)+x_{i}^{-}(t)\omega_{i}^{+}(t)).
	\end{split}
\end{gather}
The relation between the analog output of ALNN and the equilibrium of the above ODEs is $y=y^+(\infty)-y^-(\infty)$. Since $y$ is independent of the initial weights, usually generated in random, in ALNN, it is expected that the equilibrium $(y^+(\infty),y^-(\infty))^\top$ in the ODEs (\ref{eq:13}) is globally asymptotically stable, i.e., not depending on the initial condition $(y^+(0), y^-(0))^\top$. The following proposition ensures this point, which amounts to reproducing the analogous stability result given by Moorman, Samaniego, et al. (2019).
\begin{prop}
For any initial condition $(y^+(0), y^-(0))^\top$, any positive initial concentration of the catalyst species and $k, k_{\infty}>0$, the MAS (\ref{eq:6})-(\ref{eq:12}) governed by Eqs. (\ref{eq:13}) and (\ref{eq:00}) admits a unique exponentially stable equilibrium. Moreover, if $k_{\infty} \rightarrow \infty$, $y^+(t), y^-(t)$ behave as a ramp-like function at equlibrium.  
\end{prop}
\begin{pf}
The proof is very similar to that in Moorman, Samaniego, et al. (2019). $\hfill \Box$
\end{pf}

Utilizing $k=1$, when $k_{\infty} \rightarrow \infty$ and further $1/k_{\infty} \rightarrow 0$, we can obtain the equilibrium to be 
\begin{equation}
	y^{+}(\infty) =
	\begin{cases}
		p^{+}(0)-p^{-}(0), & \text{if}\quad p^{+}(0)-p^{-}(0)>0;\\
		0, & \text{if}\quad p^{+}(0)-p^{-}(0)\leqslant 0;
	\end{cases}  
\label{eq:15}
\end{equation}
\begin{equation}
	y^{-}(\infty) =
	\begin{cases}
		0, & \text{if}\quad p^{+}(0)-p^{-}(0)>0;\\
		p^{-}(0)-p^{+}(0), & \text{if}\quad p^{+}(0)-p^{-}(0)\leqslant 0.
	\end{cases}  
	\label{eq:16}
\end{equation}
Therefore, the analog output $y$ could be expressed as the steady state of the output species:
\begin{equation}
    y =
	\begin{cases}
		 y^{+}(\infty), & \text{if}\quad y>0;\\
		-y^{-}(\infty), & \text{if}\quad y<0.
	\end{cases}
\label{eq:17}
\end{equation}
This means that the MAS (\ref{eq:6})-(\ref{eq:12}) implements the forward propagation process of ALNN. For convenience, we refer to the reactions (\ref{eq:6})-(\ref{eq:12}) as the \textit{feedforward CRN} (fCRN) in the context.

\subsection{Programming Backpropagation}
The backpropagation process of ALNN includes learning, judgment, and loop, where the learning procedure means the weight update, and is our main concern here. From the weight update formula (\ref{eq:4}), we denote the error by
\begin{equation}\label{eq:gao1}
e_j(m)=d_j-y_j(m),
\end{equation} 
and also adopt the dual rail encoding technique \citep{vasic2020deep} to rewrite them (for the $i$th input, the $j$th sample and the $m$th iteration) as
\begin{equation}\label{eq:19}
\Delta \omega_{ij}^{+}(m) - \Delta \omega_{ij}^{-}(m)  = \eta (x_{ij}^{+}-x_{ij}^{-})(e_j^{+}(m)-e_j^{-}(m))
\end{equation} 
and 
\begin{equation}\label{eq:gao2}
	e_j^{+}(m)-e_j^{-}(m)=(d_j^+-d_j^-)-(y_j^+(m)-y_j^-(m)),
\end{equation} 
respectively. We further denote
\begin{align}
   e_j^{+}(m) &=d_j^++y_j^-(m), \label{eq:gao3}\\
e_j^{-}(m) & =d_j^-+y_j^+(m), \label{eq:gao4}\\
\Delta \omega_{ij}^{+}(m)&= \eta [x_{ij}^{+}e_j^{+}(m)+x_{ij}^{-}e_j^{-}(m)], \label{eq:gao5}\\
\Delta \omega_{ij}^{-}(m)  &= \eta [x_{ij}^{+}e_j^{-}(m)+x_{ij}^{-}e_j^{+}(m)]. \label{eq:gao6}
\end{align}
Hence, we need to additionally introduce the \textit{weight change species} pair $\{\Delta W_{i}^{+},\Delta W_{i}^{-}\}$, the \textit{desired output species pair} $\{D^+,D^-\}$ and the \textit{error species pair} $\{E^+,E^-\}$ for the design of reactions, which are referred to as the \textit{online learning CRN} (olCRN) in the following.

Since the weight update and the output computation take place alternatively, it needs to design the olCRN with great caution, especially for those shared species also emerging in the fCRN. Meanwhile, the parallel occurrence of all reactions also needs to be handled carefully. To address these issues, we, inspired by \cite{blount2017feedforward}, introduce two cell-like compartments as containers, where the reactions in the fCRN and olCRN take place, respectively. They are connected through so-called permeation walls, like channels, which can selectively prohibit some species (such as some macromolecular speices), release some species (like some small molecular species), or change a species on one side to a different species on the other side through reactions. The first two kinds of permeation are called physical permeation, and the third one is called chemical permeation. Utilizing this technique, it is reasonable to assume that after the fCRN reach the equilibrium, the input species, weight species,
and the output species in the feedforward compartment will react with the permeation walls to produce the corresponding auxiliary input species $\tilde{X}_i^{+}, \tilde{X}_i^{-}$, auxiliary weight species $\tilde{W}_i^{+},\tilde{W}_i^{-}$ and auxiliary output species $\tilde{Y}^{+},\tilde{Y}^{-}$ in the online learning compartment, respectively. The involved reactions, written in a simplified form, are
\begin{gather*}
	X_i^+ {\longrightarrow} \tilde{X}_i^{+},~~~ 
	W_i^+ {\longrightarrow} \tilde{W}_i^{+},~~~
	Y^+ {\longrightarrow} \tilde{Y}^{+},\\
	X_i^- {\longrightarrow} \tilde{X}_i^{-},~~~ 
	W_i^- {\longrightarrow} \tilde{W}_i^{-},~~~
	Y^- {\longrightarrow} \tilde{Y}^{-}.
\end{gather*}
where $i=1,\ldots,n$. Here and below, we ignore the identification of the number of samples. The update of samples only requires to re-assign values to $D^+$ or $D^-$ and $Y^+$ or $Y^-$ each time. After this group of permeation reactions, the auxiliary species will inherit the concentration of the corresponding original species while the original species will vanish. Note that the priority of the permeation reactions is much lower than the main reactions (\ref{eq:6})-(\ref{eq:12}) in the feedforward compartment and reactions in the online learning compartment. Namely, they will have no effect on the dynamics of the fCRN and olCRN.


We then can design the olCRN for programming the weight update. From Eq. (\ref{eq:19}), the weight update needs to individually compute the error and the weight change. We thus build the following MAS to program the error computation firstly,
\begin{align}
  	D^{+} & \xrightarrow{1} D^{+}+E^{+}, \label{eq:cata_w_1}\\
	\tilde{Y}^{+} &\xrightarrow{1} \tilde{Y}^{+}+E^{-}, \label{eq:cata_w_2}\\
	D^{-} &\xrightarrow{1} D^{-}+E^{-}, \label{eq:cata_w_3}\\
	\tilde{Y}^{-} &\xrightarrow{1} \tilde{Y}^{-}+E^{+}, \label{eq:cata_w_4}\\
	E^{+} &\xrightarrow{1}\varnothing, \label{eq:deca_e_1}\\
	E^{-} &\xrightarrow{1} \varnothing,\label{eq:deca_e_2}\\
	  D^{+}+\tilde{Y}^{+} & \xrightarrow{k_{l\infty}} \varnothing, \label{eq:ann_l_1}\\
	D^{-}+\tilde{Y}^{-} & \xrightarrow{k_{l\infty}} \varnothing, \label{eq:ann_l_2}
\end{align}
where $k_{l\infty} \gg 1$. The reactions (\ref{eq:cata_w_1}), (\ref{eq:cata_w_4}) and (\ref{eq:deca_e_1}) serves to support Eq. (\ref{eq:gao3}) while the reactions (\ref{eq:cata_w_2}), (\ref{eq:cata_w_3}) and (\ref{eq:deca_e_2}) are used to support Eq. (\ref{eq:gao4}). The last two reactions (\ref{eq:ann_l_1}) and (\ref{eq:ann_l_2}) are crucial for realizing the automatic computation of the error. We anticipate this function through a remark.

\begin{remark}
The output species enters into the online learning compartment from the feedforward compartment either in the form of $\tilde{Y}^+$ or of $\tilde{Y}^-$, so there are four combinations among the output species and the desired output species, i.e., $(D^+,\tilde{Y}^+)$, $(D^-,\tilde{Y}^-)$, $(D^+,\tilde{Y}^-)$ and $(D^-,\tilde{Y}^+)$. For the first one, if $d^+>\tilde{y}^+$, then the reactions (\ref{eq:ann_l_1}), (\ref{eq:cata_w_1}) and (\ref{eq:deca_e_1}) take place, which realizes the automatic computation of $e^+=d^+-\tilde{y}^+$; if $d^+<\tilde{y}^+$, then the reactions (\ref{eq:ann_l_1}), (\ref{eq:cata_w_2}) and (\ref{eq:deca_e_2}) take place, which realizes the automatic computation of $e^-=\tilde{y}^+-d^+$; if $d^+=\tilde{y}^+$, the backpropagation is not needed. The similar analysis may be applied to the second combination $(D^-,\tilde{Y}^-)$. For the third one, the reactions (\ref{eq:cata_w_1}) (\ref{eq:cata_w_4}) and (\ref{eq:deca_e_1}) take place, which directly performs the automatic computation of $e^+=d^++\tilde{y}^-$. The fourth one leads to $e^-=d^-+\tilde{y}^+$. 
\end{remark}

We secondly design the weight update network based on Eqs. (\ref{eq:gao5}) and (\ref{eq:gao6}), which takes
\begin{align}
    E^{+}+\tilde{X}_{i}^{+} & \xrightarrow{\eta} E^{+}+\tilde{X}_{i}^{+}+\Delta W_{i}^{+}, \label{eq:cata_update_1}\\
	E^{+}+\tilde{X}_{i}^{-} & \xrightarrow{\eta} E^{+}+\tilde{X}_{i}^{-}+\Delta W_{i}^{-}, \label{eq:cata_update_2}\\
	E^{-}+\tilde{X}_{i}^{-} & \xrightarrow{\eta} E^{-}+\tilde{X}_{i}^{-}+\Delta W_{i}^{+}, \label{eq:cata_update_3}\\
	E^{-}+\tilde{X}_{i}^{+} & \xrightarrow{\eta} E^{-}+\tilde{X}_{i}^{+}+\Delta W_{i}^{-}, \label{eq:cata_update_4}\\
	\Delta W_{i}^{+} & \xrightarrow{1} W_{i}^{+}, \label{eq:update_5}\\
	\Delta W_{i}^{-} & \xrightarrow{1} W_{i}^{-}, \label{eq:update_6}\\
	\tilde{W}_i^{+}  & \xrightarrow{1} W_{i}^{+}+\tilde{W}_i^{+}, \label{eq:update_7}\\
	\tilde{W}_i^{-} & \xrightarrow{1} W_{i}^{-}+\tilde{W}_i^{-}, \label{eq:update_8}\\
	W_{i}^{+}+ \mathcal{B} & \xrightarrow{k_1} \mathcal{B}, \label{eq:update_9}\\
	W_{i}^{-}+ \mathcal{B} & \xrightarrow{k_1} \mathcal{B}, \label{eq:update_10}
	\end{align}
where the rate constants of the first four reactions are designed to be the learning rate $\eta \in \left(0,1\right]$, and of the last two decay reactions satisfy $k_1 \beta(0)=1$. Clearly, the first four reactions are directly generated from Eqs. (\ref{eq:gao5}) and (\ref{eq:gao6}); the reaction (\ref{eq:update_5})-(\ref{eq:update_8}) are used to perform weight update, and simultaneously, to keep the original weights and balance the weight change species; the last two reactions serves to balance the weight species $W_i^+$ and $W^-_i$. Here, we introduce a new auxiliary macromolecular catalyst $\mathcal{B}$ to construct two catalytic reactions instead of the direct hydrolysis of $W_i^+$ and $W^-_i$. The main reason is that these two species also emerge in the feedforward compartment, where no hydrolysis reactions of $W_i^+$ and $W^-_i$ happen. To keep logic consistency, we construct two catalytic reactions instead in the online learning compartment.

Since the updated weights require to return to the feedforward compartment for the next iteration, we further apply the permeation technique \citep{blount2017feedforward} and the following assumption. Namely, after the olCRN reach the equilibrium, the auxiliary input species $\tilde{X}_{i}^+$,$\tilde{X}_{i}^-$ change to the original input species $X_{i}^+$,$X_{i}^-$ through chemical permeation, i.e.,   
\begin{gather*}
    \tilde{X}_i^{+} {\longrightarrow} X_i^+,~~~
    \tilde{X}_i^{-} {\longrightarrow} X_i^-;
\end{gather*}
the weight species in the olCRN return to the feedforward compartment through physical permeation (this is not contradictory to the chemical permeation for the weight species from the feedforward compartment to the online learning compartment due to asymmetry of the cell wall.); the catalyst $\mathcal{B}$ always stay in the online learning compartment; and other species, regarded as wastes, should be released to the environment. Up to now, we finish the network design to program the whole backpropagation process.

We further check the correctness of the networks from the viewpoint of dynamics. The fast reactions (\ref{eq:ann_l_1})-(\ref{eq:ann_l_2}) may induce the dynamic equations to be
\begin{equation}\label{eq:fast_system}
	\begin{split}
		\frac{d\tilde{y}^{+}(t)}{dt} & =-k_{l\infty}d^{+}(t)\tilde{y}^{+}(t)=\frac{dd^{+}(t)}{dt},\\
 		\frac{d\tilde{y}^{-}(t)}{dt} &=-k_{l\infty}d^{-}(t)\tilde{y}^{-}(t) =\frac{dd^{-}(t)}{dt}
 	\end{split}
 \end{equation}
while the slow reactions (\ref{eq:cata_w_1})-(\ref{eq:deca_e_2}) and (\ref{eq:cata_update_1})-(\ref{eq:update_10}) render  
\begin{equation}
    \begin{split}\label{eq:slow_system}
        \frac{de^{+}(t)}{dt} &=d^{+}(t)+\tilde{y}^{-}(t)-e^{+}(t),\\
		\frac{de^{-}(t)}{dt} &=d^{-}(t)+\tilde{y}^{+}(t)-e^{-}(t),\\
		\frac{d\Delta \omega_i^+ (t)}{dt} &=\eta e^+(t) \tilde{x}_i^+(t) +\eta e^-(t) \tilde{x}_i^-(t)-\Delta \omega_i^+(t),\\
		\frac{d\Delta \omega_i^- (t)}{dt} &=\eta e^+(t) \tilde{x}_i^-(t) +\eta e^-(t) \tilde{x}_i^+(t)-\Delta \omega_i^-(t),\\
		\frac{d\omega_i^+(t)}{dt} &= \Delta \omega_i^+ (t) +\tilde{\omega}_i^{+}(t)-k_1\beta(t)\omega_i^+(t),\\
		\frac{d\omega_i^-(t)}{dt} &= \Delta \omega_i^- (t) +\tilde{\omega}_i^{-}(t)-k_1\beta(t)\omega_i^-(t).
 	\end{split}
 \end{equation}
It is not difficult to check the equilibrium of the ODEs (\ref{eq:slow_system}) is expected as suggested by Eqs. (\ref{eq:gao3})-(\ref{eq:gao6}). Moreover, the last two equations together with the facts that $\mathcal{B}$ is a catalyst so $\beta(t)=\beta(0)=\beta(\infty)$ and $k_1\beta(0)=1$, render the weight update formulas  
 \begin{equation}
    \begin{split}\label{eq:weight-update}
       		\omega_i^+(t)&= \Delta \omega_i^+ (t) +\tilde{\omega}_i^{+}(t),\\
		\omega_i^-(t) &= \Delta \omega_i^- (t) +\tilde{\omega}_i^{-}(t).
 	\end{split}
 \end{equation}
 The ODEs (\ref{eq:fast_system}) serves to provide the boundary equilibrium about  $d^+,\tilde{y}^+,d^-,\tilde{y}^-$. 
 
In the following, we exhibit the designed olCRN is exponentially convergent when the concentrations of the desired output species and the auxiliary output species are given. 
\begin{theorem}
Consider the mass-action olCRN (\ref{eq:cata_w_1})-(\ref{eq:update_10}) described by Eqs. (\ref{eq:fast_system})-(\ref{eq:slow_system}). For any given positive initial condition, any given positive initial concentration of the auxiliary species and $\eta, k_1>0$, if $k_{l\infty} \rightarrow \infty$, then the solution of the ODEs (\ref{eq:fast_system})-(\ref{eq:slow_system}) exponentially converges to the boundary equilibrium of the fast system (\ref{eq:fast_system}) and the solution of the slow system (\ref{eq:slow_system}). Further, if $k_1 \beta(0)=1$, then the above olCRN can implement the update of weight $w_{ij}(m)$ by training the $j$th sample $(\bm{x_{j}},d_{j})$ at the $m$th iteration. 
\end{theorem}
\begin{pf}
See the proof in the \textbf{Appendix A}. $\hfill \Box$
\end{pf}

 


\section{Case study and Discussions}
In this section, we use a simple example to illustrate our results for implementing ALNN to fit a straight line. All codes are compiled and performed in the Matlab environment (Matlab2019a). Then, we make some discussions about the significance, challenges and points of future research on the current work. 

\begin{figure*}
  	 \begin{subfigure}{.5\textwidth}
  		\centering
  		\includegraphics[scale=0.22]{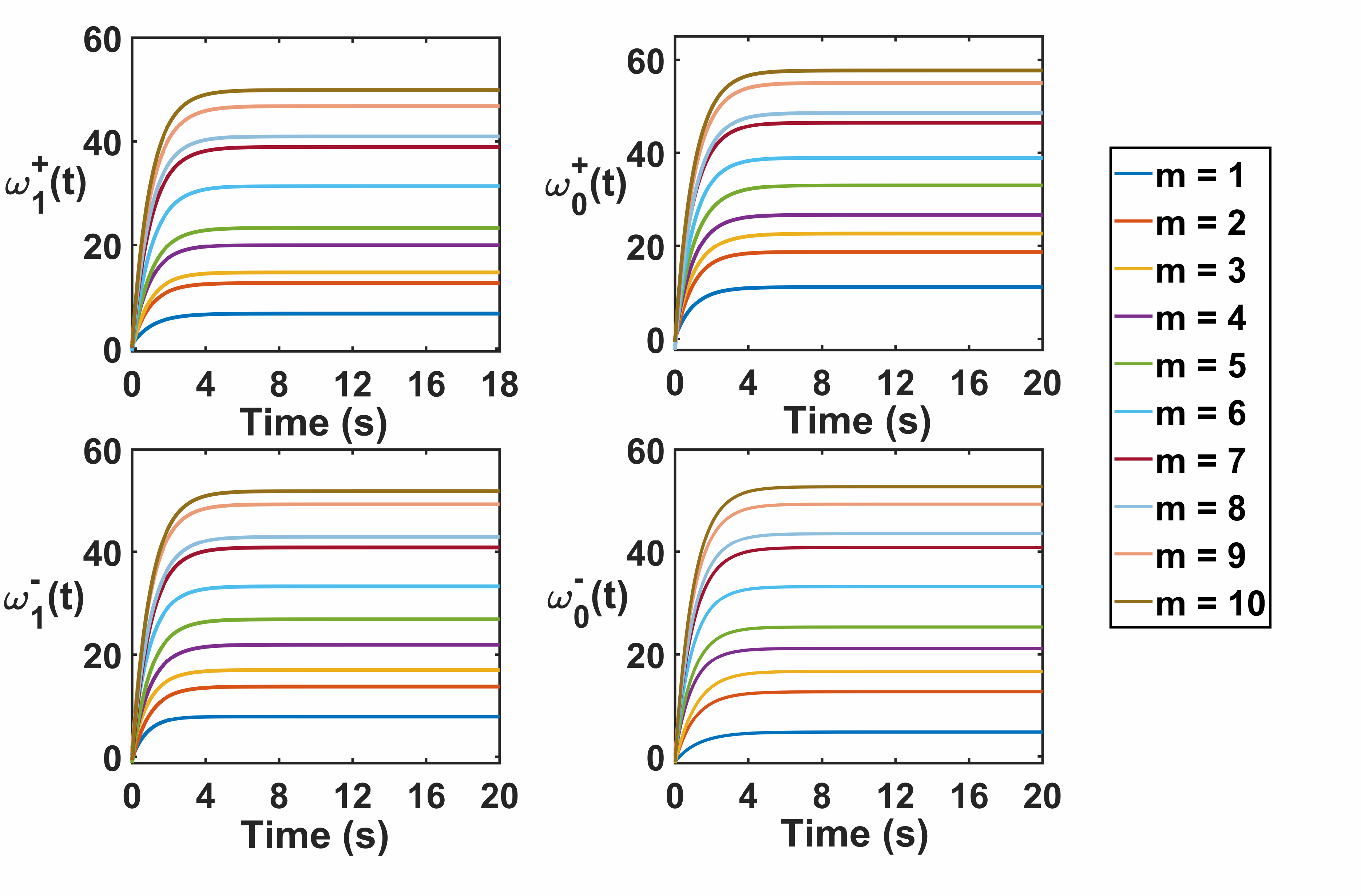}
  	 	\caption{convergence of the concentrations of weight
  	 	species}\label{fig:1a}
  	 \end{subfigure}
     \begin{subfigure}{.5\textwidth}
  	    \centering
  	    \includegraphics[width=0.75\textwidth]{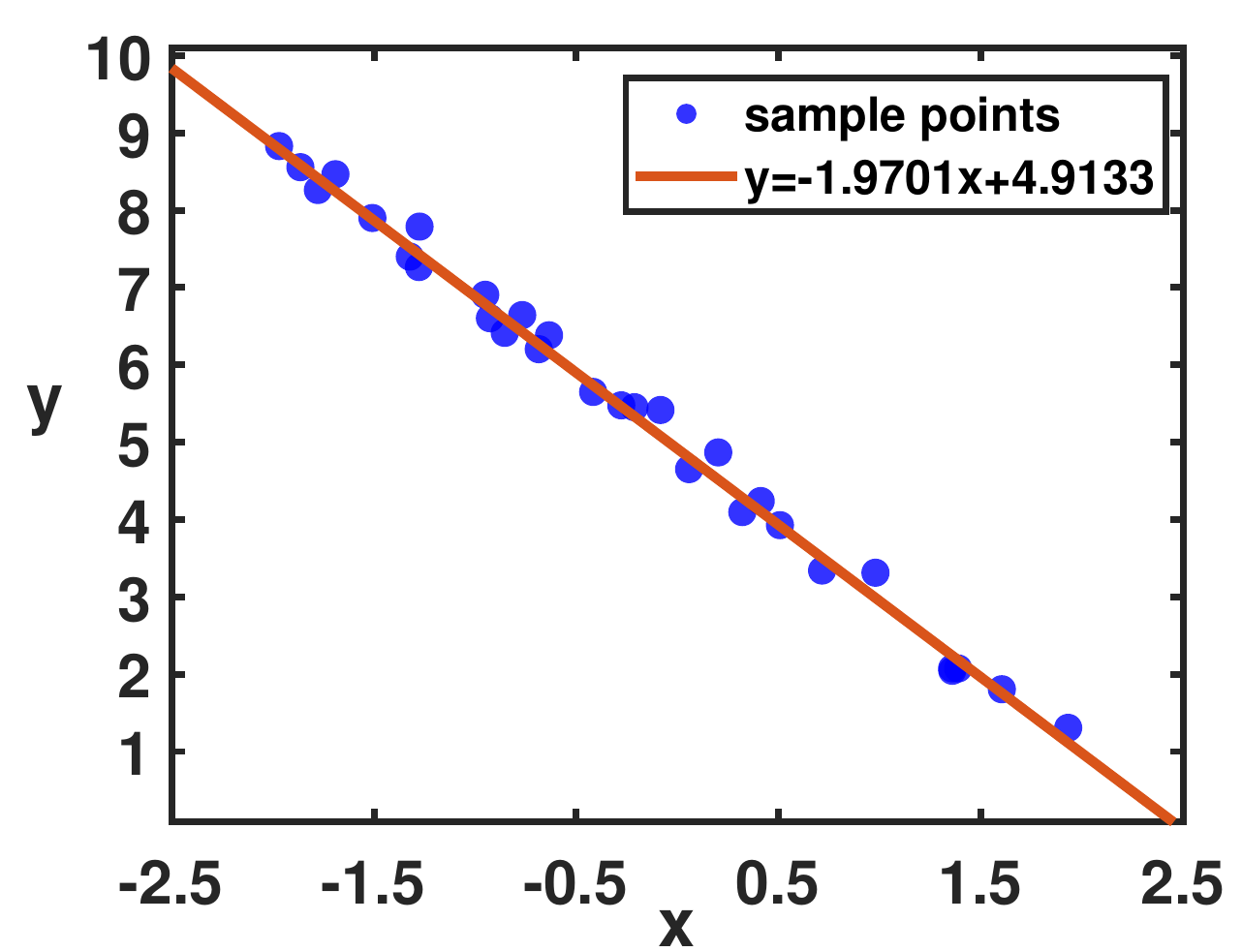}
  	    \caption{the fitting effect of samples points}\label{fig:1b}
     \end{subfigure}
     \caption{\textbf{The simulation results of biochemical ALNNs for Example 1}}
     \label{fig.1}
\end{figure*}


\begin{example}
The task is to fit the straight line $y=-2x+5$ using the CRN language. There will include one input species $X_1$, two weight species $W_0,W_1$, one output species $Y$, and one desired output species $D$. We firstly generate $30$ samples in random within $x\in [-2,2]$ bounded by error $\pm 0.3$. Let $\eta=0.2$, $k_{\infty}=k_{l\infty}=500$, the reaction time be $20$s, and the iteration termination condition be $|e|<0.1$. The positive initial value of every related species is generated following a standard normal distribution. We present the iteration process of the weight concentrations in Fig. \ref{fig:1a}, where each line looks like an exponential function implying exponential convergence. Every point in the line represents the corresponding weight concentration at different times. The final weight vector is obtained to be $\bm{\omega}=(-1.9701,4.9133)^{\top}$ at the $10$th iteration. Further, the fitting effect is displayed in Fig. \ref{fig:1b} with the straight line $y=-1.9701x+4.9133$ versus the $30$ sample points. Clearly, the CRN succeeds in implementing the function of ALNN in this linear fitting problem. 
\end{example} 

Naturally, the above example is quite simple, and ALNN is also the most basic machine learning algorithm. However, this does not prevent CRNs from implementing more complex computations and more complex algorithms. A direct extension is to implementing nonlinear neural networks, even deep learning algorithms, in which multiple layer structure, gradient computation, and other activation functions, such as Sigmoid function, ReLU function, etc., should be taken into account. And the judgment module also needs to be considered, such as designing CRNs to determine whether the iteration satisfies the termination condition while updating the weights, which is not covered in this paper. Also, one might note that in the design of the loop module, we borrow the idea of cell-like compartments from \cite{blount2017feedforward}. This technique is no doubt practicable, but it needs special biological materials as cell walls and strict experimental conditions to meet the desired requirement. Actually, the key to implementing the loop operation is to control the occurrence sequence of reactions since chemical reactions occur in parallel. A feasible solution to address this issue is to design chemical oscillators that are potential to control the occurrence time and extent of reactions. In addition, in the current reactions design, the fast reactions are often used, such as reactions (\ref{eq:12}), (\ref{eq:ann_l_1}) and (\ref{eq:ann_l_2}). The accurate analysis depends on the rate constants to be infinite. However, it is not possible in practice. How to address it is also a difficult issue. Along with this, there also exists the issue of errors when approximating the ALNN results computed in finite time by the CRNs results calculated in infinite time. All of these issues are extremely challenging, since more enormous networks will be encountered, and more masterly design and rigorous analysis on CRNs will be required.


Despite full of difficulties and challenges, there is still of great significance to develop the technique to program machine learning algorithms using CRNs, which, in our opinion, may promote the development of the following three fields: \textbf{(1) Synthetic biology}. Synthetic biology is to design, transform and even resynthesize organisms in a targeted manner through engineering design concepts, thereby creating ``synthetic organisms" endowed with unnatural functions. It usually requires assembling many modules that have some computation function like assembling an electric circuit. However, it is impossible to directly embed a certain machine learning algorithm (that can perform some computation) into the above modules. This work may provide a way for cells to use these ``algorithms" to achieve corresponding functions, such as classification, etc, in the living organism. \textbf{(2) Molecular computer}. Molecular computer is a generic term for any computational scheme which uses molecules as the carrier to accomplish computational missions instead of traditional computers. Depending on the computational and information processing capabilities of CRNs, they could be seen as another programming language CRN++ with typical branching constructing such as if/else, loops, etc \citep{vasic2020crn++}. For the completeness of algorithms, our future implementation of neural networks will include the assignment, loop, and judgment modules that are compiled to CRNs, which contributes to providing more molecular programming paradigms. \textbf{(3) Artificial intelligence (AI)}. Today, the AI techniques are developing very rapidly, whose applications are nearly covered throughout all fields. However, their effectiveness is being limited by the machine’s inability to explain its decisions and actions to users, i.e., lack of explainability and working as a black box very often. The most important reason is that there is no strong theoretical support for AI. Most of the related techniques originate from human intuition or directly imitate living behaviors. Programming machine learning algorithms using CRNs amounts to modeling the former using the latter's theory. Since the CRN theory is very rigorous, mainly being polynomial system theory that falls into the category of dynamical system, it is possible to construct the theoretical framework of AI based on the polynomial system theory.


\bibliography{ifacconf}             

     

\appendix
\section{ Proof of Theorem 1.}
\begin{pf}
All differential equations according to the online learning CRNs are written as Eqs. (\ref{eq:fast_system}), (\ref{eq:slow_system}) in the text. We interpret $\tilde{y}^+, \tilde{y}^-, d^{+}, d^{-}$ as fast variables and $e^{+}, e^{-}$ $ \Delta \omega_i^+, \Delta \omega_i^-$,
$\omega_i^+, \omega_i^-$ as slow variables in the sense of Tikhonov’s theorem (\citep{verhulst2005methods}).
Note that slow variables are regarded as constants and the slow system (\ref{eq:slow_system}) without multiplied by $k_{l\infty}$ are removed in the fast system (\ref{eq:fast_system}). Considering our participating reactions, if $\tilde{Y}$ matches $D$, fast reactions will occur and the fast system reach equilibrium immediately when $k_{l\infty} \rightarrow \infty$ according to Tikhonov’s theorem. Taking the existence of both $\tilde{Y}^+$ and $D^+$ as an example, the steady states of fast system could be written as 
	\begin{gather}
		(d^+(\infty), \tilde{Y}^+(\infty))= (d^+(0)- \tilde{Y}^+(0), 0) 
		\label{eq:supp_1}
	\end{gather}
	or
	\begin{gather}
		(d^+(\infty), \tilde{Y}^+(\infty))= (0, \tilde{Y}^+(0)-d^+(0)) 
		\label{eq:supp_2}
	\end{gather}
depending on the magnitude of the initial concentration of $\tilde{Y}$ and $D$. Thus the solution to (\ref{eq:fast_system}),(\ref{eq:slow_system}) converges to (\ref{eq:supp_1}) or (\ref{eq:supp_2}) and the solution to (\ref{eq:slow_system}).
	
Let $\bm{z_i (t)}=\bm{(z_i^{\textit{fast}}(t), z_i^{\textit{slow}}(t))}$ denotes the vector composed of all variables present in (\ref{eq:fast_system}),(\ref{eq:slow_system}). Now given positive initial concentration vector $\bm{z_i(0)}, \tilde{x}_i^+(t), \tilde{x}_i^+(t)$, $\tilde{\omega_i}^+(t),
\tilde{\omega_i}^-(t), \beta(t)$ are constants during the whole process. The slow system (\ref{eq:slow_system}) could be rewritten as a constant coefficient linear system with offset 
	\begin{equation}
		\dot{\bm{z_i^{\textit{slow}}(t)}} = U\bm{z_i^{\textit{slow}}}(t) + \bm{b} 
		\label{eq:supp_3}
	\end{equation}
where $U$ is a lower triangular matrix with diagonal elements of $-1$ and $b$ is a constant vector after substituting the equilibrium of the fast system. The equilibrium point of (\ref{eq:supp_3}) can be shifted to the origin via a change of variables with $U$ unchanged, and $U$ only admits negative eigenvalue $-1$, which ensures that the system (\ref{eq:slow_system}) has a unique globally asymptotically stable equilibrium point due to the linear system theory. So we could describe the steady states of the olCRN system below.
\begin{equation}\label{eq:supp_4}
	\begin{split}
			e^+(\infty)& =d^+(\infty)+\tilde{y}^-(\infty) \\
			e^-(\infty)& =d^-(\infty)+\tilde{y}^+(\infty) \\
			\Delta \omega_i^{+}(\infty) & = \eta e^+(\infty) \tilde{x}_i^+(0) +\eta e^-(\infty) \tilde{x}_i^-(0)\\
			\Delta \omega_i^{-}(\infty) & = \eta e^-(\infty) \tilde{x}_i^+(0) +\eta e^+(\infty) \tilde{x}_i^-(0)\\
			\omega_i^+(\infty)& = \Delta \omega_i^+ (\infty) +\tilde{\omega}_i^+(0)\\
			\omega_i^-(\infty)& = \Delta \omega_i^- (\infty) +\tilde{\omega}_i^+(0)
	\end{split}
\end{equation}
where we let $\tilde{\omega}_i^+(0), \tilde{\omega}_i^-(0)$ be $\omega_i^{+}, \omega_i^{-}$ of the previous update, respectively, according to the chemical permeation of weight species, and $\tilde{x}_i^+(0), \tilde{x}_i^-(0)$ depend on $\bm{x_j}$. Consequently, combining (\ref{eq:supp_1}), (\ref{eq:supp_2}) and (\ref{eq:supp_4}), the updated value of $\omega_i$ for the $j_{th}$ training sample is equal to the difference between $\omega_i^+(\infty)$ and $\omega_i^-(\infty)$, and the equilibrium results (\ref{eq:supp_4}) are consistent with the Eqs.(\ref{eq:gao3})-(\ref{eq:gao6}) and (\ref{eq:weight-update}) in the text. $\hfill\Box$
\end{pf}
\section{MATLAB Codes}
\subsection*{Main Program}
\begin{lstlisting}
clear all;
clc;
n = 0.2; k=500; k1=500;
load('C:\Users\Administrator\crn\data.mat')
one([1:30],1)=1; X=[rand_x,one];
row_number=size(X,1);column_number=size(X,2);
zerx=zeros(row_number,column_number);
Xp=max(X,zerx); Xn=abs(min(X,zerx));
zer=zeros(1,column_number);
w=randn(1,column_number);
wp=max(w,zer); wn=abs(min(w,zer));
W=[];E=[];
for i=1:10
    for j=1:row_number
         xp=Xp(j,:);
         xn=Xn(j,:);
         d=D(j);
     [t1,y1,yp,yn]=NNFeedforward(xp,xn,k,wp,wn,column_number);
         p=yp-yn;
     [t2,y2,e,wp,wn]=hybrid_parameter_update(xp,xn,wp,wn,p,d,n,k1,column_number);
         w=wp-wn;
         W=[W;w]; el(j)=e;
      end
         E=[E;el']; 
     if abs(e)<0.1
         break
  end
end
     \end{lstlisting}
\subsection*{Main Functions}
     \begin{lstlisting}
function [t,y1,yp,yn]= NNFeedforward(xp,xn,k,wp,wn,sn)
     X0=randn(2,1);t0=0;tfinal=20;
     [t,y1]=ode45(@(t,y) comput_output(t,y,k,sn,wp,wn,xp,xn),[t0,tfinal],X0);
     yp=y1(end,1);
     yn=y1(end,2);
end

function [t,y,e,wp_new,wn_new] = hybrid_parameter_update(xp,xn,wp,wn,p,d,n,k1,column_number)
     t0=0;tfinal=20;
     w0=randn(2+column_number*4,1);wp_new=[];wn_new=[];
end
if (d>0)&&(p>0)
     Y0=[d;p;w0];
     [t,y]=ode45(@(t,y) match_ode(t,y,wp,wn,xp,xn,n,k1,column_number),[t0,tfinal],Y0);  
     e=y(end,3)-y(end,4);
         for i=1:column_number
          wp_new=[wp_new,y(end,4+2*column_number+i)];
          wn_new=[wn_new,y(end,4+3*column_number+i)];
         end   
     else if (d<0)&&(p<0)
      Y0=[-d;-p;w0];
      [t,y]=ode45(@(t,y) match_ode2(t,y,wp,wn,xp,xn,n,k1,column_number),[t0,tfinal],Y0);
      e=y(end,3)-y(end,4);
         for i=1:column_number
          wp_new=[wp_new,y(end,4+2*column_number+i)];
          wn_new=[wn_new,y(end,4+3*column_number+i)];
         end 
     else 
          Y0=w0;
          [t,y]=ode45(@(t,y) not_match_ode(t,y,wp,wn,xp,xn,n,column_number,d,p),[t0,tfinal],Y0);
          e=y(end,1)-y(end,2);
         for i=1:column_number
          wp_new=[wp_new,y(end,2+2*column_number+i)];
          wn_new=[wn_new,y(end,2+3*column_number+i)];
         end
     end
 end
     \end{lstlisting}

\subsection*{Functions in ODE45}
\begin{lstlisting}
function dydt = comput_output(t,y,k,sn,wp,wn,xp,xn)
	sump=0;sumn=0;
  for i=1:sn
	p=wp(i)*xp(i)+wn(i)*xn(i);
	n=wp(i)*xn(i)+wn(i)*xp(i);
	sump=sump+p;
	sumn=sumn+n;
  end
	dydt=[sump-y(1)-k*y(1)*y(2);sumn-y(2)-k*y(1)*y(2)];
end
\end{lstlisting}

\begin{lstlisting}
function dydt = match_ode(t,y,vp,vn,xp,xn,n,k1,column_number)
	y_delta_wp=[];y_delta_wn=[];y_wp=[];y_wn=[];
	y_error=[y(3);y(4)];
	y1dot=-k1*y(1)*y(2);
	y2dot=-k1*y(1)*y(2);
	y_error_dot=-eye(2)*y_error+[1,0,0,1;0,1,1,0]*[y(1);0;y(2);0];
  for i=1:column_number
	y_delta_wp=[y_delta_wp;y(4+i)];
	y_delta_wn=[y_delta_wn;y(4+column_number+i)];
	y_wp=[y_wp;y(4+2*column_number+i)];
	y_wn=[y_wn;y(4+3*column_number+i)];
  end
	y_delta_wpdot=-eye(column_number)*y_delta_wp+n*[xp',xn']*y_error;
	y_delta_wndot=-eye(column_number)*y_delta_wn+n*[xn',xp']*y_error;
	y_wpdot=-eye(column_number)*y_wp+vp'+y_delta_wp;
	y_wndot=-eye(column_number)*y_wn+vn'+y_delta_wn;
	dydt=[y1dot;y2dot;y_error_dot;y_delta_wpdot;y_delta_wndot;y_wpdot;y_wndot];
end
\end{lstlisting}
    
\begin{lstlisting}
function dydt = match_ode2(t,y,vp,vn,xp,xn,n,k1,column_number)
	y_delta_wp=[];y_delta_wn=[];y_wp=[];y_wn=[];
	y_error=[y(3);y(4)];
	y1dot=-k1*y(1)*y(2);
	y2dot=-k1*y(1)*y(2);
	y_error_dot=-eye(2)*y_error+[1,0,0,1;0,1,1,0]*[0;y(1);0;y(2)];
  for i=1:column_number
	y_delta_wp=[y_delta_wp;y(4+i)];
	y_delta_wn=[ y_delta_wn;y(4+column_number+i)];
	y_wp=[y_wp;y(4+2*column_number+i)];
	y_wn=[y_wn;y(4+3*column_number+i)];
  end
	y_delta_wpdot=-eye(column_number)*y_delta_wp+n*[xp',xn']*y_error;
	y_delta_wndot=-eye(column_number)*y_delta_wn+n*[xn',xp']*y_error;
	y_wpdot=-eye(column_number)*y_wp+vp'+y_delta_wp;
	y_wndot=-eye(column_number)*y_wn+vn'+y_delta_wn;
	dydt=[y1dot;y2dot;y_error_dot;y_delta_wpdot;y_delta_wndot;y_wpdot;y_wndot];
end
\end{lstlisting}

\begin{lstlisting}
function dydt = not_match_ode(t,y,vp,vn,xp,xn,n,column_number,d,p)
  if (d>0)&&(p<0)
	dp=d;dn=0;yn=-p;yp=0;
	else
	dn=-d;dp=0;yn=0;yp=p;
  end
	y_delta_wp=[];y_delta_wn=[];y_wp=[];y_wn=[];
	y_error=[y(1);y(2)];
	y_error_dot=-eye(2)*y_error+[1,0,0,1;0,1,1,0]*[dp;dn;yp;yn];
  for i=1:column_number
	y_delta_wp=[y_delta_wp;y(2+i)];
	y_delta_wn=[ y_delta_wn;y(2+column_number+i)];
	y_wp=[y_wp;y(2+2*column_number+i)];
	y_wn=[y_wn;y(2+3*column_number+i)];
  end
	y_delta_wpdot=-eye(column_number)*y_delta_wp+n*[xp',xn']*y_error;
	y_delta_wndot=-eye(column_number)*y_delta_wn+n*[xn',xp']*y_error;
	y_wpdot=-eye(column_number)*y_wp+vp'+y_delta_wp;
	y_wndot=-eye(column_number)*y_wn+vn'+y_delta_wn;
	dydt=[y_error_dot;y_delta_wpdot;y_delta_wndot;y_wpdot;y_wndot];
end
\end{lstlisting}
\end{document}